\documentclass[11pt]{article}
\usepackage[usenames]{color}
\usepackage[colorlinks=true,
linkcolor=webgreen,
filecolor=webbrown,
citecolor=webgreen]{hyperref}

\definecolor{webgreen}{rgb}{0,.5,0}
\definecolor{webbrown}{rgb}{.6,0,0}

\usepackage{amssymb,psfig,epsfig,latexsym,graphicx,here}

\newtheorem{theorem}{Theorem}

\def\binom#1#2{{#1}\choose{#2}}

\newcommand{\eqn}[1]{(\ref{#1})}

\newcommand{\beql}[1]{\begin{equation}\label{#1}}
\newcommand{\eeq}{\end{equation}}

\makeatletter
\def\@sect#1#2#3#4#5#6[#7]#8{\ifnum #2>\c@secnumdepth
     \def\@svsec{}\else
     \refstepcounter{#1}\edef\@svsec{\csname the#1\endcsname.\hskip .75em }\fi
     \@tempskipa #5\relax
      \ifdim \@tempskipa>\z@
        \begingroup #6\relax
          \@hangfrom{\hskip #3\relax\@svsec}{\interlinepenalty \@M #8\par}%
        \endgroup
       \csname #1mark\endcsname{#7}\addcontentsline
         {toc}{#1}{\ifnum #2>\c@secnumdepth \else
                      \protect\numberline{\csname the#1\endcsname}\fi
                    #7}\else
        \def\@svsechd{#6\hskip #3\@svsec #8\csname #1mark\endcsname
                      {#7}\addcontentsline
                           {toc}{#1}{\ifnum #2>\c@secnumdepth \else
                             \protect\numberline{\csname the#1\endcsname}\fi
                       #7}}\fi
     \@xsect{#5}}
\def\@begintheorem#1#2{\it \trivlist \item[\hskip \labelsep{\bf #1\ #2.}]}

\def\section{\@startsection {section}{1}{\z@}{-3.5ex plus -1ex minus
 -.2ex}{2.3ex plus .2ex}{\normalsize\bf}}
\makeatother
\makeatletter
\def\subsection{\@startsection {subsection}{1}{\z@}{-3.5ex plus -1ex minus
 -.2ex}{2.3ex plus .2ex}{\normalsize\bf}}

\makeatother

\begin{document}
\begin{center}

{\large\bf Acyclic Digraphs and Eigenvalues of $(0,1)$--Matrices} \\
\vspace*{+.2in}
\end{center}

\begin{center}
Brendan D. McKay, 
Department of Computer Science, 
Australian National University, \\ 
Canberra, ACT 0200, AUSTRALIA  
\smallskip

Fr\'ed\'erique~E.~Oggier\footnote{This 
work was carried out during F.~E.~Oggier's 
visit to AT\&T Shannon Labs during the summer 
of 2003. She thanks the Fonds National Suisse, 
Bourses et Programmes d'\'Echange for support.}, 
D\'epartement de Math\'ematiques, \\ 
Ecole Polytechnique F\'ed\'erale de Lausanne, 
1015 Lausanne, SWITZERLAND 
\smallskip

Gordon F. Royle, 
Department of Computer Science \& Software Engineering, 
University of \\
Western Australia, 
35 Stirling Highway, 
Crawley, WA 6009, AUSTRALIA  
\smallskip

N.~J.~A.~Sloane\footnote{To whom
correspondence should be addressed. [Email: njas@research.att.com, 
phone: 973 360 8415, fax: 973 360 8178.]}, 
Internet and Network Systems Research Department, 
AT\&T Shannon Labs, \\ 
180 Park Avenue, 
Florham Park, NJ 07932--0971, USA 
\smallskip

Ian M. Wanless, 
Department of Computer Science, 
Australian National University, \\ 
Canberra, ACT 0200, AUSTRALIA 
\smallskip

Herbert S. Wilf, 
Mathematics Department, 
University of Pennsylvania, \\ 
Philadelphia, PA 19104--6395, USA  

\bigskip
October 24, 2003 \bigskip

\end{center}

\begin{center}
{\bf Abstract}
\end{center}

We show that the number of acyclic directed graphs with $n$
labeled vertices is equal to the number of $n\times n$ $(0, 1)$--matrices 
whose eigenvalues are positive real numbers.

\vspace{0.8\baselineskip}

\vspace{0.8\baselineskip}

\section{Weisstein's conjecture}
A calculation was recently made by Eric W. Weisstein 
of Wolfram Research, Inc.,
to count the real $n\times n$ matrices of $0$'s and $1$'s all
of whose eigenvalues are real and positive.
The resulting sequence of values, viz.,
\[1,3,25,543,29281\]
(for $n = 1, 2, \ldots, 5$)
was then observed to coincide with the beginning
of sequence 
\htmladdnormallink{A003024}{http://www.research.att.com/cgi-bin/access.cgi/as/njas/sequences/eisA.cgi?Anum=A003024}
in \cite{OEIS},
which counts acyclic digraphs with $n$ labeled vertices.
Weisstein conjectured that the sequences were in fact identical,
and we prove this here.

\noindent Notation. A ``digraph'' means a graph with at most one edge
directed from vertex $i$ to vertex $j$, for $1 \le i \le n, 1 \le j \le n$.
Loops and cycles of length two are permitted, but parallel edges are
forbidden. ``Acyclic'' means there are no cycles of any length.
\begin{theorem}
For each $n=1,2,3,\dots$, the number of acyclic directed graphs with
$n$ labeled vertices is equal to the number of
$n\times n$ matrices of $0$'s and $1$'s whose eigenvalues are
positive real numbers.
\end{theorem}

\paragraph{Proof.}
Suppose we are given an acyclic directed graph $G$. 
Let $A=A(G)$ be its vertex adjacency matrix. 
Then $A$ has only $0$'s on the diagonal,
else cycles of length $1$ would be present. 
So define $B=I+A$, and note that $B$ is also a matrix 
of $0$'s and $1$'s. We claim $B$ has only positive eigenvalues.

Indeed, the eigenvalues will not change if we renumber the
vertices of the graph $G$ consistently with the partial order
that it generates. But then $A=A(G)$ would be strictly upper 
triangular, and $B$ would be upper triangular with $1$'s on
the diagonal. Hence all of its eigenvalues are equal to $1$.

Conversely, let $B$ be a $(0,1)$--matrix whose
eigenvalues are all positive real numbers. Then we have
\begin{eqnarray}\label{Eq1}
1&\ge& \frac{1}{n}\mathrm{Trace}(B)\qquad\quad\ (\mathrm{since\ all\ }B_{i,i}\le 1) \nonumber \\
&=&\frac{1}{n}(\lambda_1+\lambda_2+\dots+\lambda_n) \nonumber \\
&\ge&\left(\lambda_1\lambda_2\dots\lambda_n\right)^{\frac{1}{n}}\qquad
(\mbox{by~the~arithmetic-geometric~mean~inequality}) \nonumber \\
&=&\left(\det{B}\right)^{\frac{1}{n}} \nonumber \\
&\ge&1\qquad\qquad\quad\qquad\quad (\mathrm{since}\ \det{B}\ \mathrm{is\ a\ positive\ integer}).
\end{eqnarray}
Since the arithmetic and geometric means of the eigenvalues are equal,
the eigenvalues are all equal, and in fact all $\lambda_i(B)=1$.

Now regard $B$ as the adjacency matrix of a digraph $H$,
which has a loop at each vertex. Since
\[\mathrm{Trace}(B^k)=\sum_{i=1}^n\lambda_i^k=\sum_{i=1}^n1=n,\]
for all $k$, the number of closed walks in $H$, of each length $k$, is $n$.

Since the trace of $B$ is equal to $n$, all diagonal entries
of $B$ are $1$'s. Thus we account for all $n$ of the
closed walks of length $k$ that exist in the graph $H$
by the loops at each vertex.
There are no closed walks of any length that use an edge of $H$
other than the loops at the vertices.

Put $A=B-I$. Then $A$ is a $(0,1)$--matrix that is the
adjacency matrix of an acyclic digraph. $\Box$

Remark. The only related result we have found in the literature
is the theorem \cite[p. 81]{CDS95}
that a digraph G contains no cycle if and only if all eigenvalues
of the adjacency matrix are 0.

\section{Corollaries.}
The proof also establishes the following results.

(i) Let $B$ be a $(0,1)$--matrix whose eigenvalues are
all positive real numbers. Then the eigenvalues are in fact all equal to $1$.
The only symmetric $(0,1)$--matrix with positive
eigenvalues is the identity.

(ii) Let $B$ be an $n\times n$ matrix with integer entries
and Trace$(B)\le n$. Then $B$ has all eigenvalues real
and positive if and only if $B=I+N$, where $N$ is nilpotent.

(iii) If a digraph contains a cycle,
then its adjacency matrix has an eigenvalue which is
zero, negative, or strictly complex.
In fact, a more detailed argument, not given here,
shows that if the length of the shortest cycle is at least $3$, then
there is a strictly complex eigenvalue.

(iv) The eigenvalues of a digraph consist of
$n-k$ $0$'s and  $k$ $1$'s if and only if 
the digraph is acyclic apart from $k$ loops.

(v) Define two matrices $B_1$, $B_2$ to be {\em equivalent} if
there is a permutation matrix $P$ such that
$P' B_1 P ~=~ B_2$.
Then the number of equivalence classes of $n \times n$ (0,1)--matrices
with all eigenvalues positive is equal to the number
of acyclic digraphs with $n$ unlabeled vertices.
(These numbers form sequence
\htmladdnormallink{A003087}{http://www.research.att.com/cgi-bin/access.cgi/as/njas/sequences/eisA.cgi?Anum=A003087} in \cite{OEIS}.) \\
Proof. Two labeled graphs $G_1$, $G_2$
with adjacency matrices $A(G_1)$, $A(G_2)$ correspond
to the same unlabeled graph if and only if
there is a permutation matrix $P$ such that
$P' A(G_1) P ~=~ A(G_2)$. The result now follows
immediately from the theorem. $\Box$

(vi) Let $B$ be an $n \times n$ $(-1, +1)$--matrix
with all eigenvalues real and positive. Then $n=1$ and $B=[1]$. \\
Proof. The argument that led to \eqn{Eq1} still applies
and shows that all the eigenvalues are $1$, $\det B = 1$
and Trace$(B) = n$. By adding or subtracting the first row
of $B$ from all other rows we can clear the first column,
obtaining a matrix
$$
C ~=~
\left[ \begin{array}{cc} 1 & \ast \\ \bf{0} & D \end{array} \right] \,,
$$
where $\bf{0}$ is a column of $0$'s and $D$ is an $n-1 \times n-1$ matrix
with entries $-2, 0, +2$ and $\det D = \det C = \det B = 1$. Hence $2^{n-1}$
divides $1$, so $n=1$. $\Box$

It would be interesting to investigate the 
connections between matrices and graphs in other cases--for example if
the eigenvalues are required only to be real and nonnegative
(see sequences
\htmladdnormallink{A086510}{http://www.research.att.com/cgi-bin/access.cgi/as/njas/sequences/eisA.cgi?Anum=A086510},
\htmladdnormallink{A087488}{http://www.research.att.com/cgi-bin/access.cgi/as/njas/sequences/eisA.cgi?Anum=A087488} in \cite{OEIS} for the initial values),
or if the entries are $-1$, $0$ or $1$
(\htmladdnormallink{A085506}{http://www.research.att.com/cgi-bin/access.cgi/as/njas/sequences/eisA.cgi?Anum=A085506}).

\section{Bibliographic remarks}
Acyclic digraphs were first counted by Robinson \cite{Rob70, Rob73}, and
independently by Stanley \cite{Sta73}: if $R_n$ is the
number of acyclic digraphs with $n$ labeled vertices, then
$$
R_n ~=~ \sum_{k=1}^{n} (-1)^{k+1} {\binom{n}{k}} 2^{k(n-k)} R_{n-k} ~,
$$
for $n \ge 1$, with $R_0 = 1$, and
\[\sum_{n=0}^{\infty}R_n\frac{x^n}{2^{{n\choose
2}}n!}=\left[\sum_{n=0}^{\infty}(-1)^n\frac{x^n}{2^{{n\choose
2}}n!}\right]^{-1} ~.\]
The asymptotic behavior is
\[R_n ~\sim~ n!\frac{2^{{n\choose 2}}}{Mp^n} ~,\]
where $p=1.488\ldots$ and $M=0.474\ldots$.

The asymptotic behavior of $R(n,q)$, the number of these graphs that have $q$
edges, was found by Bender \textit{et al.} \cite{BRRW, BR88}, and the number
that have specified numbers of sources and sinks has been found by Gessel
\cite{Ges96}.

\end{document}